\documentclass[10pt,reqno]{amsart}
\usepackage{amsmath}
\usepackage{amsthm}
\usepackage{amssymb}
\usepackage{amsfonts}
\usepackage{latexsym}
\usepackage{setspace}

\newcommand{\C}{ \mathbb{C}}
\newcommand{\D}{ \mathbb{D}}

\newcommand{\norm}[1]{\| #1 \|}
\newcommand{\inner}[1]{\langle #1 \rangle}
\newcommand{\del}{\partial}

\newcommand{\E}{\mathcal{E}}
\newcommand{\F}{\mathcal{F}}

\newcommand{\W}{\mathcal{W}}

\newcommand{\h}{\mathcal{H}}

\newcommand{\minimatrix}[4]{\begin{pmatrix} #1 & #2 \\ #3 & #4 \end{pmatrix}  }
\newcommand{\megamatrix}[9]{\begin{pmatrix} #1 & #2 & #3 \\ #4 & #5 & #6 \\ #7 & #8 & #9\end{pmatrix}  }

\newcommand{\rank}{\operatorname{rank}}
\renewcommand{\vec}[1]{#1}

\renewcommand{\phi}{\varphi}
\newcommand{\model}{H^2 \ominus \phi H^2}

\newtheorem{Corollary}{Corollary}
\newtheorem{Theorem}{Theorem}

\theoremstyle{definition}
\newtheorem*{Definition}{Definition}

\newtheorem{Example}{Example}

\allowdisplaybreaks
\begin{document}
    \title{Some New Classes of Complex Symmetric Operators}

    \author{Stephan Ramon Garcia}
    \address{   Department of Mathematics\\
            Pomona College\\
            Claremont, California\\
            91711 \\ USA}
    \email{Stephan.Garcia@pomona.edu}
    \urladdr{http://pages.pomona.edu/\textasciitilde sg064747}
    
    \author{Warren R. Wogen}
    \address{Department of Mathematics\\
	CB \#3250, Phillips Hall\\
	Chapel Hill, NC 27599}
    \email{wrw@email.unc.edu}
    \urladdr{http://www.math.unc.edu/Faculty/wrw}

    
    \thanks{First author partially supported by National Science Foundation Grant DMS-0638789.}

    \begin{abstract}
    	We say that an operator $T \in B(\h)$ is complex symmetric if there
	exists a conjugate-linear, isometric involution $C:\h\rightarrow\h$ 
	so that $T = CT^*C$.  We prove that 
	binormal operators, operators that are algebraic of degree two (including all idempotents),	
	and large classes of rank-one perturbations of normal operators are complex symmetric.
	From an abstract viewpoint, these results explain why the compressed shift
	and Volterra integration operator are complex symmetric.
	Finally, we attempt to describe all complex symmetric partial isometries,	
	obtaining the sharpest possible statement given only the data $(\dim \ker T, \dim \ker T^*)$.
    \end{abstract}

\maketitle

\section{Introduction}

Throughout this note, $\h$ will denote a separable complex Hilbert space and all
operators considered will be bounded.  We first require a few preliminary definitions:

\begin{Definition}
	A \emph{conjugation} is a conjugate-linear operator $C:\h \rightarrow \h$, 
	that is both \emph{involutive} ($C^2 = I$) and \emph{isometric}.
	We say that a bounded linear operator $T \in B(\h)$ is \emph{$C$-symmetric}
	if $T = CT^*C$ and \emph{complex symmetric} if there exists a conjugation $C$ with respect to which $T$
	is $C$-symmetric.
\end{Definition}

It is not hard to see that $T$ is a complex symmetric operator 
if and only if $T$ is unitarily equivalent to a symmetric matrix with complex entries,
regarded as an operator acting on an $l^2$-space of the appropriate dimension (see 
\cite[Sect. 2.4]{GarciaCCO} or \cite[Prop. 2]{GarciaPutinar}).

The class of complex symmetric operators includes all normal operators, 
operators defined by Hankel matrices, compressed Toeplitz operators 
(including finite Toeplitz matrices and the compressed shift), 
and the Volterra integration operator.  We refer the reader to \cite{GarciaPutinar,GarciaPutinar2} (or \cite{GarciaCCO}
for a more expository pace) for further details.  
Other recent articles concerning complex symmetric operators include \cite{Chevrot, Gilbreath}.

In this note, we exhibit several additional classes of complex symmetric operators.  In particular, we
establish that\medskip
\begin{enumerate}\addtolength{\itemsep}{0.5\baselineskip}
	\item All binormal operators are complex symmetric (Theorem \ref{TheoremBinormal})
		and that $n$-normal operators that are not complex symmetric exist for each $n \geq 3$
		(Example \ref{ExampleNilpotent}).
		
	\item Operators that are algebraic of degree two are complex symmetric (Theorem \ref{TheoremAlgebraic}).
		This includes all idempotents and all operators that are nilpotent of order $2$.
		
	\item Large classes of rank-one perturbations of normal operators are complex symmetric
		(Theorem \ref{TheoremNormal}).
		On abstract grounds, this explains why the compressed
		shift operator (Example \ref{ExampleShift})
		and Volterra integration operator (Example \ref{ExampleVolterra})
		are complex symmetric.
		
	\item We attempt to describe all complex symmetric partial isometries,	
		obtaining the sharpest possible statement (Theorem \ref{TheoremPartial})
		given only the data $(\dim \ker T, \dim \ker T^*)$.
\end{enumerate}

\section{Binormal Operators and $n$-normal Operators}

\begin{Definition}
	An operator $T \in B(\h)$ is called \emph{binormal} if $T$ is unitarily equivalent to an operator of the form
	\begin{equation}\label{eq-Block}
		\minimatrix{N_{11}}{N_{12}}{N_{21}}{N_{22}}
	\end{equation}
	where the entries $N_{ij}$ are commuting normal operators.  
	More generally, we say that $T$ is \emph{$n$-normal} 
	if $T$ is unitarily equivalent an $n \times n$ 
	operator matrix whose entries are commuting
	normal operators.
\end{Definition}

Needless to say, each $n \times n$ scalar matrix trivially defines an $n$-normal operator
on $\C^n$.  For further information concering binormal and $n$-normal operators, we refer the reader to \cite{HL,RR}.

The main theorem of this section is:

\begin{Theorem}\label{TheoremBinormal}
	If $T \in B(\h)$ is a binormal operator, then 
	$T$ is a complex symmetric operator.  This result is sharp in the sense that 
	if $n \geq 3$, then there exists an $n$-normal operator
	that is not a complex symmetric operator.
\end{Theorem}

\begin{proof}
	We focus our attention on the first statement since the second
	will follow from the construction of explicit examples
	(see Example \ref{ExampleNilpotent}).
	Given an operator of the form \eqref{eq-Block}, 
	the Spectral Theorem asserts that we may assume that
	each $N_{ij}$ is a multiplication operator $M_{u_{ij}}$ on a Lebesgue space $L^2(\mu)$ where  
	$\mu$ is a Borel measure on $\C$ with compact support $\Delta$ and that the corresponding
	symbols $u_{ij}$ belong to $L^{\infty}(\mu)$.  
	To simplify our notation, we will henceforth identify multiplication 
	operators $M_u$ with their symbols $u$.
	
	Without loss of generality, we may further restrict our attention to operators 
	on $L^2(\mu)^{(2)}$ (the two-fold inflation of $L^2(\mu)$) of the form
	\begin{equation}\label{eq-Triangular}
		T = \minimatrix{ u_{11} }{ u_{12} }{0}{u_{22}}
	\end{equation}		
	since any binormal operator is unitarily equivalent 
	to an operator of form \eqref{eq-Triangular} \cite[Thm.~7.20]{RR}.
%	since Schur's Theorem on unitary upper-triangularization ensures that
%	any binormal operator is unitarily equivalent an operator of the form 
	
	Let us denote by $E$ the subset of $\Delta$ upon which $u_{11} = u_{22}$:
	\begin{equation*}
		E = \{ \,z \in \Delta\, :\, u_{11}(z) = u_{22}(z)\, \,\text{$\mu$-a.e.}\, \}.
	\end{equation*}
	Letting $\chi_E$ denote the characteristic function of $E$, we note that 
	the subspace $\E_1 = \chi_E L^2(\mu)$ and its orthogonal complement
	$\E_2 = \chi_{\Delta \backslash E} L^2(\mu)$ are both
	reducing subspaces for $M_z:L^2(\mu)\rightarrow L^2(\mu)$, the operator
	of multiplication by the independent variable.
	In particular, their inflations $\E_1^{(2)}$ and $\E_2^{(2)}$ are both reducing subspaces for $T$ and
	we see that
	\begin{equation*}
		T = T | _{\E_1^{(2)}} \oplus T | _{\E_2^{(2)}}.
	\end{equation*}
	Since the direct sum of complex symmetric operators is complex symmetric,
	we need only consider the following two special cases:	\medskip
	\begin{enumerate}\addtolength{\itemsep}{0.5\baselineskip}
		\item $u_{11} = u_{22}$ $\mu$-a.e.
		\item $u_{11} \neq  u_{22}$ $\mu$-a.e.
	\end{enumerate}
	
	\bigskip
	\noindent \textsc{Case} (i):  Suppose that $u_{11} = u_{22}$
	$\mu$-a.e.  In this case, we may write \eqref{eq-Triangular} as
	\begin{equation*}
		\minimatrix{u}{v}{0}{u}
	\end{equation*}
	where $u,v \in L^{\infty}(\mu)$.  One can immediately verify that
	$T$ is $C$-symmetric with respect to the conjugation
	$C(f_1,f_2) = (\overline{f_2}, \overline{f_1})$ on $L^2(\mu)^{(2)}$.
	
	\bigskip
	\noindent \textsc{Case} (ii):  Suppose that $u_{11} \neq  u_{22}$ $\mu$-a.e.
	In this case, $T$ has the form
	\begin{equation}\label{eq-BlockUV}
		 \minimatrix{u_1}{v}{0}{u_2}
	\end{equation}
	where $u_1 \neq u_2$ $\mu$-a.e.  	
	Let $F$ denote the subset of $\Delta$ upon which $v$ vanishes
	and observe that 
	$T = T | _{\F_1^{(2)}} \oplus T | _{\F_2^{(2)}}$
	where $\F_1 = \chi_F L^2(\mu)$ and $\F_2 = \chi_{\Delta \backslash F} L^2(\mu)$.
	Since $v$ vanishes on $F$, it follows from \eqref{eq-BlockUV}
	that $T | _{\F_1^{(2)}}$ is normal and hence
	complex symmetric.  On the other hand, 
	$T | _{\F_2^{(2)}}$ is an operator of the form \eqref{eq-BlockUV} where $v$
	is $\mu$-a.e. nonvanishing.
	Without loss of generality, we may therefore assume that
	$v$ does not vanish on a set of positive $\mu$-measure.  	
	
	Since $u_1 - u_2$ and $v$ are nonvanishing $\mu$-a.e., we may define
	a unimodular function $\gamma$ by the formula
	\begin{equation}
		\gamma = \frac{v}{|v|} \cdot \frac{ |u_1 - u_2|}{(u_1 - u_2)}. \label{eq-Gamma}
	\end{equation}
	Letting
	\begin{align}
		a(z) &= \frac{ \gamma |u_1 - u_2| }{ \sqrt{   |u_1 - u_2|^2 +|v|^2} }, &
		b(z) &= \frac{ |v| }{ \sqrt{  |u_1 - u_2|^2 + |v|^2} },\label{eq-ab}
	\end{align}
	we note that the operator
	\begin{equation*}
		U = \minimatrix{  a}{b}{b}{ - \overline{a} }
	\end{equation*}
	on $L^2(\mu)^{(2)}$ is unitary since $b$ is real and $|a|^2 + |b|^2 = 1$  $\mu$-a.e.
	
	Let $Jf = \overline{f}$ denote the canonical conjugation on $L^2(\mu)$ and let $K = J^{(2)}$
	denote its two-fold inflation:
	\begin{equation*}
		K = \minimatrix{J}{0}{0}{J}.
	\end{equation*}
	Clearly $K$ is a conjugation on $L^2(\mu)^{(2)}$ and
	a short computation shows that $U^* = KUK$ (i.e. $U$ is a $K$-symmetric operator).
	
	We now claim that $C = UK$ is a conjugation on $L^2(\mu)^{(2)}$.  
	Since $C$ is obviously conjugate-linear
	and isometric, we need only verify that $C^2 = I$:
	\begin{equation*}
		C^2 = (UK)(UK) = U(KUK) = U U^* = I.
	\end{equation*}
	Thus $C$ is a conjugation on $L^2(\mu)^{(2)}$, as claimed.

	We conclude the proof by showing that $T$ is $C$-symmetric.  We will do this by
	directly verifying that $CT^* = TC$.  First note that
	\begin{align}
		TC 
		&= TUK  \nonumber\\
		&= \minimatrix{u_1}{v}{0}{u_2} \minimatrix{  a}{b}{b}{ - \overline{a} } K\nonumber\\
		&= \minimatrix{a u_1 + bv}{b u_1  - \overline{a}v}{b u_2 }{-\overline{a}u_2} K. \label{eq-Compare1}
	\end{align}
	On the other hand, we also have
	\begin{align}
		CT^*
		&= UKT^*\nonumber\\
		&=  \minimatrix{  a}{b}{b}{ - \overline{a} } K \minimatrix{ \overline{u_1} }{0}{ \overline{v} }{ \overline{u_2} }\nonumber\\
		&=  \minimatrix{  a}{b}{b}{ - \overline{a} } \minimatrix{ u_1 }{0}{ v }{ u_2 } K\nonumber\\
		&= \minimatrix{a u_1 + bv}{  bu_2  }{b u_1  - \overline{a}v}{-\overline{a}u_2} K. \label{eq-Compare2}
	\end{align}
	To verify the equality of \eqref{eq-Compare1} and \eqref{eq-Compare2}, 
	we need only show that $bu_2 = bu_1 - \overline{a} v$.  However, the
	preceding equation follows directly from \eqref{eq-Gamma} and \eqref{eq-ab}.	
\end{proof}

One might regard Theorem \ref{TheoremBinormal} as a generalization
of the following well-known result (alternate proofs of which can be found in 
\cite[Cor.~3.3]{Chevrot}, \cite[Ex.~6]{GarciaPutinar}, or \cite[Cor.~3]{Tener}):

\begin{Corollary}\label{Corollary2x2}
	Every linear operator on $\C^2$ is complex symmetric.  In other words, every 
	$2 \times 2$ matrix is unitarily equivalent to a symmetric matrix with complex entries.
\end{Corollary}

In order to verify the second claim of Theorem \ref{TheoremBinormal},
we must exhibit examples of $n$-normal operators ($n \geq 3$)
that are not complex symmetric.  
The following example does just this.

\begin{Example}\label{ExampleNilpotent}
	We first claim that the operator $T:\C^3 \rightarrow \C^3$ defined by the matrix
	\begin{equation}\label{eq-3x3}
		\megamatrix{0}{a}{0}{0}{0}{b}{0}{0}{0}
	\end{equation}
	(with respect to the standard basis)
	is a complex symmetric operator if and only if $ab = 0$ or $|a| = |b|$.  
	There are several possible case to investigate:
	(i) $a=0$ or $b = 0$,
	(ii) $|a| = |b| \neq 0$,
	(iii) $a \neq 0$, $b \neq 0$, and $|a| \neq |b|$.
	In particular, the final case yields $3$-normal operators that are not complex symmetric.
	
	\medskip
	
	\noindent \textsc{Case} (i):  If $a = 0$ or $b = 0$, 
	then $T$ is the direct sum of complex symmetric operators
	by Corollary \ref{Corollary2x2}.
	
	\medskip
	
	\noindent \textsc{Case} (ii):  If $|a| = |b| \neq 0$, then a short computation shows that
	\eqref{eq-3x3} is unitarily equivalent to a constant multiple of a 
	$3 \times 3$ nilpotent Jordan matrix.  It follows from 
	\cite[Example 4]{GarciaPutinar} or \cite[Sect. 2.2]{GarciaCCO}
	that $T$ is a complex symmetric operator.
	
	\medskip
	
	\noindent \textsc{Case} (iii):  Let $a \neq 0$, $b \neq 0$, and $|a| \neq |b|$ and
	suppose toward a contradiction that $T = CT^*C$ for some conjugation $C$.
	Let $\vec{e}_1, \vec{e}_2, \vec{e}_3$ denote the standard basis for $\C^3$ and 
	observe that $\vec{e}_1$ and $\vec{e}_3$ span the eigenspaces
	of $T$ and $T^*$, respectively, 	corresponding to the eigenvalue zero.
	Since $T^i \vec{x} = \vec{0}$ if and only if $(T^*)^i (C\vec{x}) = \vec{0}$,
	we see that 
	\begin{align*}
		C \vec{e}_1 = \alpha_1 \vec{e}_3, \quad
		C \vec{e}_2 = \alpha_2 \vec{e}_2, \quad
		C \vec{e}_3 = \alpha_3 \vec{e}_1
	\end{align*}
	where $\alpha_1, \alpha_2, \alpha_3$ are certain unimodular constants.
	The desired contradiction will arise from computing $\norm{T\vec{e}_2}$ in two different ways.
	On one hand, we have
	\begin{equation*}
		\norm{T \vec{e}_2}
		= \norm{T^*C \vec{e}_2}
		= \norm{T^*(\alpha_2 \vec{e}_2)}
		= \norm{T^* \vec{e}_2}
		= \norm{ (0,0,\overline{b}) }
		= |b|.
	\end{equation*}
	On the other hand, we also have
	$\norm{T \vec{e}_2} =  \norm{ (a,0,0) } = |a|$.  However, this contradicts the fact that
	$|a| \neq |b|$.  Therefore $T$ is not a complex symmetric operator.

	\medskip

	If $n > 3$, then we can use the preceding ideas to construct examples of
	$n$-normal operators that are not complex symmetric.
	Specifically, let $T:\C^3 \rightarrow \C^3$ be defined as in \eqref{eq-3x3},
	with $ab \neq 0$ and $|a| \neq |b|$ as in Case (iii).  The operator
	$T \oplus I$ on $\C^n$, where $I$ denotes the identity operator on $\C^{n-3}$,
	is trivially $n$-normal.  The same 
	argument used in Case (iii) reveals that $T \oplus I$ is not complex symmetric.
\end{Example}

We remark that matrices of the form \eqref{eq-3x3} arose
in a related unitary equivalence problem.
Consideration of Jordan canonical forms reveals that each 
$n \times n$ matrix is \emph{similar} to its transpose.  
On the other hand, the matrix
\begin{equation*}
	\megamatrix{0}{2}{0}{0}{0}{1}{0}{0}{0}
\end{equation*}
is not \emph{unitarily equivalent} to its transpose \cite[Prob. 159]{HalmosLA}.

We close this section with a corollary:

\begin{Corollary}
	If $N$ is a normal operator having spectral multiplicity $\leq 2$ and
	if $T$ is an operator commuting with $N$, then $T$ is a 
	complex symmetric operator.
\end{Corollary}

\begin{proof}
	If $N$ is a normal operator having spectral multiplicity $\leq 2$,
	then we may write $N = N_1 \oplus N_2^{(2)}$ where $N_1$ and
	$N_2$ are mutually singular $*$-cyclic normal operators 
	\cite[Thm.~IX.10.20]{Conway}.
	Moreover, we also have
	$T = T_1 \oplus T_2$, where $T_1$ commutes 
	with $N_1$ and $T_2$ commutes with $N_2^{(2)}$
	(see the discussion following \cite[Lem.~IX.10.19]{Conway} or \cite{Gilbreath}).
	From this we immediately see that 
	$T_1$ is normal and $T_2$ is binormal \cite[Cor.~IX.6.9, Prop.~IX.6.1.b]{Conway}.  
	It then follows from Theorem \ref{TheoremBinormal} that
	$T$ is a complex symmetric operator.  
\end{proof}

In the preceding, observe that 
if $N$ has spectral multiplicity two, then 
the conjugation corresponding to the operator $T$ depends on $T$ (as well as $N$).\medskip

Our next corollary asserts that any square root (normal or otherwise) 
of a normal operator is itself a complex symmetric operator:

\begin{Corollary}
If $T^2$ is normal, then $T$ is a complex symmetric operator.
\end{Corollary}

\begin{proof}
This follows immediately from Theorem \ref{TheoremBinormal}
and the fact that $T$ must be of the form
\begin{equation*}
T = A \oplus \minimatrix{B}{C}{0}{-B}
\end{equation*}
where $A$ and $B$ are normal and $C$ is a positive operator
that commutes with $B$ \cite[Thm.~1]{RR2}.
\end{proof}

\section{Algebraic Operators}

\begin{Definition}
	An operator $T \in B(\h)$ is \emph{algebraic} if
	$p(T) = 0$ for some polynomial $p(z)$.  The \emph{degree} of an algebraic
	operator is defined to be the degree of the polynomial $p(z)$ of least degree
	for which $p(T) = 0$.	
\end{Definition}

Although the following theorem is essentially a corollary of Theorem \ref{TheoremBinormal},
we choose to state it as a theorem since it will have several useful corollaries of its own.

\begin{Theorem}\label{TheoremAlgebraic}
	If $T\in B(\h)$ is algebraic of degree $\leq 2$, then $T$ is a complex symmetric operator.  
	This result is sharp in the sense that for each finite $n \geq 3$ and for
	each $\h$ satisfying $\dim \h \geq n$, there exists an algebraic
	operator on $\h$ of degree $n$ that is not a complex symmetric operator.
\end{Theorem}

\begin{proof}
	The first statement follows from Theorem \ref{TheoremBinormal}
	and an old lemma of Brown  \cite[Lem. 7.1]{Brown} that asserts that
	if $T$ is algebraic of degree $\leq 2$, then $T$ is binormal.
	Suppose now that $3 \leq n \leq \dim \h$ and consider the operator 
	$T \oplus D$ where $T$ has a matrix of the form \eqref{eq-3x3}
	with $ab \neq 0$ and $|a| \neq |b|$ and 
	$D$ is a diagonal operator  chosen so that $T \oplus D$ is algebraic of degree $n$.
	An argument similar to that used in Case (iii) of Example \ref{ExampleNilpotent} shows
	that this operator is not complex symmetric.	
\end{proof}

Two particular classes of operators stand out for special consideration:

\begin{Corollary}\label{Cor-Idempotent}
	Let $T \in B(\h)$.  If $T$ is idempotent (i.e. $T^2 = T$) or 
	nilpotent of order $2$ (i.e. $T^2 = 0$), then $T$ is a complex symmetric operator.
\end{Corollary}

A direct proof of the second portion of Corollary \ref{Cor-Idempotent}, involving the explicit construction
of the associated conjugation, can be found in \cite{GarciaATCSO}.  
Yet another basic class of operators that happen to be complex symmetric are
the rank-one operators:

\begin{Corollary}\label{Cor-Rank}
	If $T \in B(\h)$ and $\rank(T) =1$, then $T$ is a complex symmetric operator.
\end{Corollary}

\begin{proof}
	Any rank-one operator $T$ is of the form
	$Tf = \inner{f,v}u$ for certain vectors $u,v$ (this operator
	is frequently denoted $u \otimes v$).
	Since $T^2 - \inner{u,v} T = 0$,  it follows from Theorem \ref{TheoremAlgebraic}
	that $T$ is a complex symmetric operator.
\end{proof}

It is important to note that although every operator on $\C^2$ is a complex
symmetric operator (Corollary \ref{Corollary2x2}), 
there are certainly operators having \emph{rank} two that are not complex symmetric
operators (Example \ref{ExampleNilpotent}).

\section{Perturbations of Normal Operators}

In light of Corollary \ref{Cor-Rank} and the fact that
all normal operators are complex symmetric 
(see \cite[Ex. 2.8]{GarciaCCO} or \cite[Sect. 4.1]{GarciaPutinar}),
it is natural to attempt to identify those 
rank-one perturbations of normal operators that are also complex symmetric.

\begin{Theorem}\label{TheoremNormal}
	If $N \in B(\h)$ is a normal operator, $U$ is a unitary operator in
	$\W^*(N)$ (the von Neumann algebra generated by $N$), $a \in \C$,
	and $v \in \h$, then the operator $T = N + a(Uv \otimes v)$ is a
	complex symmetric operator.
\end{Theorem}

	\begin{proof}
		We may without loss of generality
		assume that $N$ is a $*$-cyclic normal operator with cyclic vector $v$.
		Otherwise let $\h_1$ denote the reducing subspace of $N$
		generated by $v$ and let $\h_2 = \h_1^{\perp}$.  Now write
		$N = N_1 \oplus N_2$ relative to the orthogonal
		decomposition $\h = \h_1 \oplus \h_2$.  It follows that $N_1$ is
		$*$-cyclic and, since $\h_1$ reduces $U$, we have
		$T = (N_1 + a(U_1v \otimes v)) \oplus N_2$ where $U_1 = U|_{\h_1}$ belongs
		to $\W^*(N_1)$.
		
		By the Spectral Theorem,
		we may further presume that $N = M_z$, the operator
		of multiplication by the independent variable on a Lebesgue space
		$L^2(\mu)$, that $v$ is the constant function $1$, and that
		$U = M_{\theta}$, the operator of multiplication by some unimodular
		function $\theta$ in $L^{\infty}(\mu)$.	
		At this point, a straightforward computation shows that $Cf = \theta \overline{f}$
		is a conjugation on $L^2(\mu)$ with respect
		to which both $M_z$ and $\theta \otimes 1$ are $C$-symmetric.  
	\end{proof}

On an abstract level, the preceding theorem indicates that 
compressed shifts are complex symmetric operators.
In other words, starting from the fact that the Aleksandrov-Clark unitary operators are complex symmetric, 
we can directly derive the fact that the compressed shift is also complex symmetric.
In essence, this is the reverse of the path undertaken in \cite{GarciaCCO}
(which the reader may consult for further details concerning the following example).

%%%%%%%%%%%%%%%%
\begin{Example}\label{ExampleShift}
	Let $\phi$ denote a nonconstant inner function and let $H^2$ denote the 
	Hardy space on the unit disk $\D$.
	For each $\lambda$ in the open unit disk $\D$, we define the unit vectors
	\begin{align}
	  b_{\lambda}(z) &= \frac{z - \lambda}{1 - \overline{\lambda}z},\\
	  k_{\lambda}(z) &= \sqrt{ \frac{1 - |\lambda|^2}{1 - |\phi(\lambda)|^2} } \cdot \frac{1 - \overline{\phi(\lambda)}\phi(z) }{1 - \overline{\lambda}z},\\
	  q_{\lambda}(z) &= \sqrt{ \frac{1 - |\lambda|^2}{1 - |\phi(\lambda)|^2} } \cdot \frac{\phi(z) - \phi(\lambda)}{z - \lambda}.
	\end{align}
	In particular, the function $k_{\lambda}$ is a normalized reproducing kernel for the so-called
	\emph{model space} $\model$.

	For each unimodular constant $\alpha$, we define 
	the \emph{generalized Aleksandrov-Clark operator} by setting 
	\begin{equation*}
		U_{\lambda} f =
		\begin{cases}
			b_{\lambda}f      & f \perp q_{\lambda} \\
			\alpha k_{\lambda} & f = q_{\lambda}.
		\end{cases}
	\end{equation*}
	Each $U_{\lambda}$ is $C$-symmetric with respect to the conjugation
	(defined in terms of boundary functions) $[Cf](z) = \overline{fz}\phi$ on $\model$.
	Moreover, we also note that $q_{\lambda} = Ck_{\lambda}$ for each $\lambda$.

	By Theorem \ref{TheoremNormal}, it follows that the operator
	\begin{equation}\label{EquationUlambda}	
		S_{\lambda} = U_{\lambda} - (\alpha + \phi(\lambda))(k_{\lambda} \otimes q_{\lambda})
	\end{equation}
	is complex symmetric since it is of the form
	$U_{\lambda} + a(U_{\lambda}v \otimes v)$ where $a$ is a complex constant
	and $v = q_{\lambda}$.  More specifically, tracing through the proof of Theorem \ref{TheoremNormal},
	we expect that $S_{\lambda}$ will be $C$-symmetric with respect to the $C$ described above.

	The significance of this example lies in the fact that, 
	for the choice $\alpha = -\phi(\lambda) / |\phi(\lambda)|$,
	the operator \eqref{EquationUlambda} turns out to be 
	\begin{equation}\label{EquationSlambda}
		S_{\lambda}f = P_{\phi}(b_{\lambda}f),
	\end{equation}
	the compression of the operator $M_{b_{\lambda}}:H^2\rightarrow H^2$ to the subspace $\model$.	  
	Here $P_{\phi}$ denotes the orthogonal projection from $H^2$ onto $\model$.  
	The operator $S_0 f = P_{\phi}(zf)$ is commonly known as the \emph{compressed shift} 
	or \emph{Jordan model} operator corresponding to $\phi$.
	In summary, purely operator-theoretic considerations guarantee that the operators
	$S_{\lambda}$ are complex symmetric.  We refer the reader to \cite{GarciaCCO} and \cite{Sarason}
	for more information.
\end{Example}

	In fact, the preceding example can be greatly generalized (without any reference to
	function theory whatsoever).  Given a contraction $T \in B(\h)$, there is a unique
	decomposition $\h = \h_0 \oplus \h_u$ where $\h_0$ and $\h_u$
	are both $T$-invariant, $T|_{\h_u}$ is unitary, and $T|_{\h_0}$ is completely
	nonunitary (i.e., $T|_{\h_0}$ is not unitary when restricted to any of its
	invariant subspaces).
	The operator $D_T = (I - T^*T)^{1/2}$ is called the \emph{defect operator} of $T$ and
	the \emph{defect spaces} of $T$ are defined to be the subspaces $\mathcal{D}_T = \overline{ D_T \h}$ and
	$\mathcal{D}_{T^*} = \overline{ D_{T^*} \h}$.  The \emph{defect indices} of $T$
	are the numbers $\del_T = \dim \mathcal{D}_T$ and 
	$\del_{T^*} = \dim \mathcal{D}_{T^*}$.
	We say that $T \in C_{0\cdot}$ if $T^n \rightarrow 0$ (SOT) and that $T \in C_{\cdot 0}$
	if $T^* \in C_{0 \cdot}$.  Finally, we also define $C_{00} = C_{0 \cdot} \cap C_{\cdot 0}$.

	It turns out that any Hilbert space contraction with defect indices $\del_T =\del_{T^*} =1$ is complex symmetric.
	Although this is known (see \cite[Cor.~3.2]{Chevrot} for a general proof)
	and easy to prove if $T \in C_{00}$ (see \cite[Thm.~5.1]{GarciaCCO}, 
	\cite[Prop.~3]{GarciaPutinar}, \cite[Lem.~2.1]{Sarason}), 
	we are able to establish this result in the abstract -- without the use of characteristic functions and complex
	analysis.

	\begin{Corollary}
		If $T \in B(\h)$ is a contraction such that  $\del_T = \del_{T^*} =1$,
		then $T$ is a complex symmetric operator.
	\end{Corollary}

	\begin{proof}
		Since $\del_T = 1$, it follows that $I - T^*T = u \otimes u$ for some nonzero 
		vector $u$.  If $x$ is any vector orthogonal to $u$, then we have
		\begin{equation*}
			\norm{x}^2 - \norm{Tx}^2 = \inner{(I - T^*T)x,x} = |\inner{u,x}|^2 = 0.
		\end{equation*}
		Thus $T$ is isometric on a subspace of $\h$ having codimension one.
		Similarly, we see that $I - TT^*$ is also of rank-one whence
		$I - TT^* = v \otimes v$ for some nonzero vector $v$.  Putting this together, we find that	
		$T = T|_{ u^{\perp}} + c(u \otimes v)$ for some constant $c$.  In particular, there exists a unitary $U$
		such that $T = U + c'(u \otimes v)$ is a rank-one-perturbation of $U$.
		Since $T$ is of the form $T = U + a(Uv\otimes v)$ where $U$ is unitary, it follows from
		Theorem \ref{TheoremNormal} that $T$ is a complex symmetric operator.
	\end{proof}
	
	Following Theorem \ref{TheoremNormal} in another direction, we obtain the following:

\begin{Corollary}\label{Cor-AB}
	Let $T = A + iB$ denote the Cartesian decomposition of $T \in B(\h)$ (i.e. $A = A^*$ and $B = B^*$).
	If $\rank(A) =1$ or $\rank(B)=1$, then $T$ is a complex symmetric operator.
\end{Corollary}

\begin{proof}
	If $A$ has rank-one, then $A = a(v \otimes v)$ for some $a \in \mathbb{R}$
	and $v \in \h$.  Apply Theorem \ref{TheoremNormal}, with $N = iB$ and $U = I$.
\end{proof}

The preceding corollary easily furnishes many examples of non-normal complex symmetric operators.
Indeed, if $A$ is an arbitrary selfadjoint operator and $B$ is a rank-one selfadjoint operator
that does not commute with $A$, then $T = A +iB$ is a non-normal complex symmetric operator.
Despite the apparent simplicity of such a recipe, nontrivial examples abound.
Consider the following example:

\begin{Example}\label{ExampleVolterra}
	It is well-known that the \emph{Volterra integration operator} 
	\begin{equation*}
		[Vf](x) = \int_0^x f(y)\,dy
	\end{equation*}
	on $L^2[0,1]$ is a rank-one selfadjoint perturbation of a skew-selfadjoint
	operator (see \cite{GarciaAAEPRI} or \cite[Pr. 188]{Halmos} for further details).
	Indeed, a short computation shows that the selfadjoint component of $V$ is
	\begin{equation*}
		A = \tfrac{1}{2}(V+V^*) = \int_0^1 f(y)\,dy = \tfrac{1}{2}(1 \otimes 1),
	\end{equation*}
	where the $1$ above denotes the constant function. 
	By Corollary \ref{Cor-AB}, we conclude that $V$ is a complex symmetric operator.
	In fact, $V = CV^*C$ where $C$ denotes the conjugation $[Cf](x) = \overline{f(1-x)}$ 
	on $L^2[0,1]$ (see \cite[Ex. 6]{GarciaPutinar2} and
	\cite[Sect. 4.1]{GarciaCCO}).
\end{Example}

Setting $U = I$ in Theorem \ref{TheoremNormal} provides a 
generalization of Corollary \eqref{Cor-AB}:

\begin{Corollary}
	If $N \in B(\h)$ is a normal operator,
	$P$ is a rank-one orthogonal projection, 
	and $a \in \C$, 
	then $T = N + aP$ is a complex symmetric operator.
\end{Corollary}

It is important to note that not every rank-one perturbation of 
a normal operator will be complex symmetric (unless $\dim \h = 2$ -- see
Corollary \ref{Corollary2x2}).
In fact, even a rank-one perturbation of an orthogonal projection
may fail to be complex symmetric:

\begin{Example}\label{Example3x3}
	We claim that the operator $T:\C^3 \rightarrow \C^3$ defined 
	by the matrix 
	\begin{equation*}
		\megamatrix{0}{0}{1}{0}{1}{1}{0}{0}{0}
	\end{equation*}
	(with respect to the standard basis) is not a complex symmetric operator.  
	First observe that the eigenspaces of $T$ (and hence of $T^*$)
	for the eigenvalues $0$ and $1$ are both one dimensional.
	The eigenspaces of $T$ corresponding to the eigenvalues $0$ and $1$
	are spanned by the unit vectors
	$\vec{v}_0 = (1,0,0)$ and $\vec{v}_1 = (0,1,0)$, respectively.
	The eigenspaces of $T^*$ corresponding to the eigenvalues $0$ and $1$
	are spanned by the unit vectors
	$\vec{w}_0 = (0,0,1)$ and $\vec{w}_1 = (0, \frac{1}{\sqrt{2}}, \frac{1}{\sqrt{2}})$, respectively.
	If $C$ is a conjugation such that $T = CT^*C$, then 
	$0 = |\inner{ \vec{v}_0, \vec{v}_1}| 
		= | \inner{ C \vec{v}_1, C \vec{v}_0} |
		= |\inner{ \vec{w}_1, \vec{w}_0 } |
		= \tfrac{1}{\sqrt{2}}$,
	which is absurd.
\end{Example}

\section{Partial Isometries}

In this section, we attempt to classify those partial isometries that are complex symmetric.
This question is related to the preceding material in the sense that if $\phi(0) = 0$ in Example \ref{ExampleShift}, 
then the corresponding compressed shift operator is a complex symmetric partial isometry.  

Given only the dimensions of the kernels of a partial isometry and its adjoint, the 
following theorem is as definitive as possible:

\begin{Theorem}\label{TheoremPartial}
	Let $T \in B(\h)$ be a partial isometry.\smallskip
	\begin{enumerate}\addtolength{\itemsep}{0.5\baselineskip}
		\item If $\dim \ker T = \dim \ker T^* \leq 1$, 
			then $T$ is a complex symmetric operator,
		
		\item If $\dim \ker T \neq \dim \ker T^*$, 
			then $T$ is not a complex symmetric operator.
			
		\item If $2 \leq \dim \ker T = \dim \ker T^* \leq \infty$, then
			either possibility can (and does) occur.
	\end{enumerate}
\end{Theorem}

\begin{proof}
	(i) If $\dim \ker T = \dim \ker T^* = 0$, then $T^*T = TT^* = I$ whence $T$
	is unitary and hence complex symmetric.
	Suppose that $T$ is a partial isometry satisfying $\dim \ker T = \dim \ker T^* = 1$
	and that $\ker T$ and $\ker T^*$ are spanned by the unit vectors 
	$v$ and $w$, respectively.  Since
	the operator $N  = T + w \otimes v$ is unitary, it follows that
	$T = N - Nv \otimes v$ is a complex symmetric operator
	by Theorem \ref{TheoremNormal}.
	For different proof, see \cite[Cor. 3.2]{Chevrot}.
	\medskip

	\noindent (ii) We prove the contrapositive.
	If $T$ is $C$-symmetric, then it is easy to see that $Tx = 0$ if and only if $T^*(Cx) = 0$.
	Therefore $C$ furnishes an isometric, conjugate-linear bijection between $\ker T$ and $\ker T^*$
	whence $\dim \ker T = \dim \ker T^*$.
	\medskip

	\noindent (iii) This portion of the theorem follows upon consideration of several examples.
	It is trivial to produce complex symmetric partial isometries with $\dim \ker T = \dim \ker T^* = n$
	for any $n$.  In fact, $T = I \oplus 0$, where $0$ is the zero operator on an $n$-dimensional Hilbert space,
	is such an example.  	On the other hand, finding partial isometries that are not complex symmetric 
	when $2 \leq n \leq \infty$ is more involved.
	
	For the remainder of this proof,
	we choose not to distinguish between matrices and the operators
	they induce (with respect to the standard basis).
	We must first study a certain auxiliary matrix that will be used in our construction.
	Specifically, we intend to prove that 
	\begin{equation*}
		A = 	
		\begin{pmatrix}
			0 & \frac{1}{2} & 0\\
			0 & 0 & \frac{1}{4} \\
			1 & 0 & 0
		\end{pmatrix}
	\end{equation*}
	is not a complex symmetric operator.  This will follow from a careful study of the eigenstructures of $A$
	and $A^*$.  First, note that the eigenvalues of $A$ are
	\begin{equation*}	
		\lambda_1 = \tfrac{1}{2}, \qquad 
		\lambda_2 = -\tfrac{1}{4} + i \tfrac{\sqrt{3}}{4}, \qquad
		\lambda_3 = -\tfrac{1}{4} - i \tfrac{\sqrt{3}}{4},
	\end{equation*}
	and that these are also the eigenvalues of $A^*$.  
	A straightforward computation shows that corresponding unit eigenvectors of $A$ are
	\begin{align*}
		v_1 &= \tfrac{1}{\sqrt{6}} (1,\,1,\,2),\\
		v_2 &= \tfrac{1}{2\sqrt{6}} (-1+i \sqrt{3},\,-1-i \sqrt{3},\,4),\\
		v_3 &= \tfrac{1}{2 \sqrt{6}}( 1+i \sqrt{3},\,1-i \sqrt{3},\,-4).
	\end{align*}
	Since $A$ has three distinct eigenvalues, it follows that 
	$v_1,v_2,v_3$ must be sent to unimodular scalar multiples of the
	corresponding unit eigenvectors
	\begin{align*}
		w_1 &= \tfrac{1}{3}(2,\,2,\,1),\\
		w_2 &= \tfrac{1}{3}(-1-i \sqrt{3},\,-1+i \sqrt{3},\,1),\\
		w_3 &= \tfrac{1}{3}(-1+i \sqrt{3},\,-1-i \sqrt{3},\,1)
	\end{align*}
	of $A^*$.  Now observe that
	\begin{equation*}
		| \inner{v_1,v_2} | = |\inner{v_2,v_3}|  = |\inner{v_3,v_1}| = \tfrac{1}{2},
	\end{equation*}
	whereas
	\begin{equation*}
		| \inner{w_1,w_2} | = |\inner{w_2,w_3}|  = |\inner{w_3,w_1}| = \tfrac{1}{3}.
	\end{equation*}
	The same argument used in Example \ref{Example3x3} now reveals that 
	$A$ is cannot be a complex symmetric operator.
	
	We are now ready to construct our desired partial isometry.  Noting that	
	\begin{equation*}
		A^*A = \megamatrix{1}{0}{0}{0}{\frac{1}{4}}{0}{0}{0}{\frac{1}{16}},
	\end{equation*}
	we see that if $2 \leq n \leq \infty$, then 
	the $n \times 3$ matrix 
	\begin{equation*}
		B = 
		\begin{pmatrix}
			0 & \frac{ \sqrt{3}}{2} & 0\\
			0 & 0  & \frac{ \sqrt{15} }{4}\\
			0 & 0 & 0 \\
			\vdots & \vdots & \vdots\\
			0 & 0 & 0
		\end{pmatrix}
	\end{equation*}
	satisfies $A^*A + B^*B = I$ (the $3 \times 3$ identity matrix).
	The $(n+3) \times (n+3)$ matrix
	\begin{equation*}
		T = \minimatrix{A}{0}{B}{0}
	\end{equation*}
	is a partial isometry since $T^*T$ is the orthogonal projection
	\begin{equation*}
		P = \minimatrix{I}{0}{0}{0}.
	\end{equation*}
	Since it is clear from the construction of $T$ that $\dim \ker T = \dim \ker T^* = n$, we need only prove
	that $T$ is not a complex symmetric operator. 
	
	Suppose toward a contradiction that $T$ is $C$-symmetric.		
	By \cite[Thm. 2 \& Cor. 1]{GarciaPutinar2}, we may
	write $T = CJP$ where $J$ is an auxiliary conjugation that commutes with $P$.
	Since $JP = PJ$ we find that
	\begin{align*}
		J(PT)J = J(PCJP)J = PJCP = T^*P = (PT)^*
	\end{align*}
	whence $PT$ is $J$-symmetric.  However,
	\begin{equation*}
		PT =  \minimatrix{A}{0}{0}{0}
	\end{equation*}
	and the same argument that showed that $A$ was not a complex symmetric operator
	also shows that $PT$ is not a complex symmetric operator.  This contradiction shows that
	our partial isometry $T$ is not a complex symmetric operator, as desired.
\end{proof}

We remark that in the final paragraph of the proof, 
we could have appealed to the fact that the Aluthge transform of 
a complex symmetric operator is also complex symmetric \cite{GarciaATCSO}.

Based upon the preceding material, we can prove that 
every partial isometry on a three-dimensional Hilbert space is complex symmetric:

\begin{Corollary}\label{Corollary3d}
	If $\dim \h = 3$, then every partial isometry $T\in B(\h)$ is complex symmetric.
\end{Corollary}

\begin{proof}
	Suppose that $\dim \h =3$ and that $T$ is a partial isometry on $\h$.
	There are four cases to discuss:	
	\begin{enumerate}\addtolength{\itemsep}{0.5\baselineskip}
		\item If $\dim \ker T = 0$, then $T$ is unitary and thus complex symmetric.
			Indeed, the Spectral Theorem asserts that $T$ has a diagonal
			matrix representation with respect to some orthonormal basis of $\h$.
			
		\item If $\dim \ker T = 1$, then $T$ is complex symmetric by (i) of Theorem \ref{TheoremPartial}.
			The condition $\dim \ker T = \dim \ker T^*$ holds trivially since $\h$ is finite-dimensional.
		
		\item If $\dim \ker T = 2$, then $\rank(T)=1$.  By Corollary \ref{Cor-Rank}, it follows that 
			$T$ is complex symmetric.

		\item If $\dim \ker T = 3$, then $T = 0$ and the result is trivial. \qedhere
	\end{enumerate}
\end{proof}

Based on the construction used in the proof of Theorem \ref{TheoremPartial}, it is clear
that many partial isometries that are not complex symmetric exist if the dimension of the underlying
Hilbert space is $\geq 5$.
On the other hand, we were for a considerable time  
unable to determine whether all partial isometries on a four-dimensional Hilbert 
space are complex symmetric (they are).  In this setting, the method of Corollary \ref{Corollary3d}
suffices to resolve all but the case $\dim \ker T = 2$.  

Significant numerical evidence in favor of the assertion that all partial isometries
on a four-dimensional Hilbert space are complex symmetric has been produced 
by J.~Tener \cite{Tener}.  We refer the reader to
\cite{CSPI} for the resolution of this problem.

\end{document}